\documentclass[12pt]{article}

\usepackage{amsmath, amssymb, amsthm, amscd}
\usepackage{graphics}

\newtheorem{theorem}      {Theorem}

\newtheorem{proposition}  {Proposition}

\newtheorem{tauconjecture}{The $\mathbf{\tau}$ Conjecture for Polynomials}

\theoremstyle{definition}
\newtheorem{example}      {Example}

\newtheorem{definition}   {Definition}

\newcommand{\pnp}          {\ensuremath {\mathcal P \ne \mathcal {NP}\ }}
\newcommand{\upnp}         {\ensuremath {\mathcal {UP} \not \supseteq 
                             \mathcal {NP} \cap \mathcal K}\ }

\newcommand{\size}	   {\mathrm{Size}}

\newcommand{\hn}           {HN}
\newcommand{\yes}          {^{\mathrm{yes}}}
\newcommand{\no}          {^{\mathrm{no}}}
\newcommand{\nullstellensatz}
                           {\ensuremath{(\hn, \hn \yes)}}

\title{Ultimate Polynomial Time}
\author{Gregorio Malajovich\footnote{
Departamento de Matem\'atica Aplicada, Universidade Federal
do Rio de Janeiro. Caixa Postal 68530, CEP 21945, Rio de Janeiro,
RJ, Brasil. e-mail: gregorio@labma.ufrj.br.
On leave at MSRI, 1000 Centennial Drive, Berkeley CA
94720-5070. e-mail: gregorio@msri.org
}}

\date{January 21, 1999}
\begin{document}
\maketitle

\begin{abstract}
	The class $\mathcal{UP}$ of `ultimate polynomial time'
	problems over $\mathbb C$ is introduced; it contains the
	class $\mathcal P$ of polynomial time problems over
	$\mathbb C$.
\par
	The $\tau$-Conjecture for polynomials implies that $\mathcal{UP}$ 
	does not contain the class of non-deterministic polynomial 
	time problems definable without constants over $\mathbb C$. 
	This latest statement implies that $\mathcal P \ne \mathcal{NP}$
	over $\mathbb C$.
\par
	A notion of `ultimate complexity' of a problem is suggested. It
	provides lower bounds for the complexity of structured problems.
\end{abstract}

\section{Introduction}
\begin{sloppypar}
	A model of Computation and Complexity over a ring was developed
	in~\cite{BSS} and ~\cite{BCSS}, generalizing the classical
	$\mathcal{NP}$-completeness theory~\cite{GAREY-JOHNSON}. Of
	particular interest is the model of Complexity over the ring
	$\mathbb C$ of complex numbers. \end{sloppypar}
\par
	In the model of complexity over $\mathbb C$, a machine is
	allowed to input, to output and to store complex numbers, 
	to compute polynomials and to branch on equality (See
	the textbook~\cite{BCSS}
	for background). 
	This model shares some of the features of the classical 
	(Turing) model of computation (There is a discussion
	in~\cite{MM}). 
	It is known~\cite{KOIRAN97,LECTURE}
	that the hypothesis $\mathcal {BPP} \not \supseteq \mathcal{NP}$
	in the Turing setting implies $\pnp$ over $\mathbb C$. ($\mathcal 
	{BPP}$ stands for Bounded Probability Polynomial Time. If 
	$\mathcal {BPP}$ would happen to contain $\mathcal{NP}$, then
	there would be polynomial time randomized algorithms for such
	tasks as factorizing large integers or breaking most modern
	cryptographic systems).
\medskip
\par
	In \cite{Shub-Smale, BCSS, SmaleXXI}, the hypothesis $\pnp$ over
	the Complex numbers was related to a number-theoretical
	conjecture. Define a straight-line program as a list
\[
	s_0 = 1\ ,\ s_1 = x\ ,\  s_2\ , \ \cdots\ , \ s_{\tau}
\]
	where $s_i$ is, for $i\ge 2$, either $s_j+s_k$, $s_j-s_k$ or
	$s_j s_k$, for some $j, k < i$. Each $s_i$ is thus a polynomial
	in $x$. The straight-line program is said to {\em compute} the
	polynomial $s_{\tau}(x)$.
\par
	Given a polynomial $f \in \mathbb Z [x]$, the quantity $\tau(f)$
	is defined as the smallest $\tau$ such that there exists a
	straight-line program $s_0, \cdots, s_{\tau}$ computing $f(x)$. 
	For instance, $\tau(x^{2^n}-1) = 2+n$. Similarly, if $g \in
	\mathbb Z[x_1, \cdots , x_n]$, then $\tau(g)$ is the minimal
	length of a straight-line program $s_0=1, s_1=x_1, \cdots,
	s_n=x_n, s_{n+1}, \cdots, s_{\tau}=g(x)$.
\begin{tauconjecture}
	There is a constant $a > 0$ such that for any univariate polynomial
	$f \in \mathbb Z[x]$,
\[
	n(f) < \tau(f)^a
\]
	where $n(f)$ is the number of integer zeros of $f$, without 
	multiplicity.
\end{tauconjecture}
\medskip
\par
	It is known~\cite{BCSS} that the $\tau$-Conjecture for polynomials
	implies $\pnp$ over $\mathbb C$. A main step towards this result is
	the fact that, if the $\tau$-Conjecture is true, then the
	polynomials
\[
	p_d(x) = (x-1) (x-2) \cdots (x-d)
\]
	are {\em ultimately hard to compute}. This means that there cannot
	be constants $a$ and $b$ such that, for any degree $d$, for some
	non-zero polynomial $f$ (depending on $d$), we would have 
\[
	\tau \left(\ p_d (x) f(x) \ \right) < a \left( \log_2 d \right) ^b
\]
\par
	Therefore, all non-zero multiples of $p_d$ are hard to compute,
	hence the wording ultimately hard.
\medskip
\par
	The goal of this paper is to define a new complexity class
	$\mathcal{UP}$, of {\em ultimate polynomial time} problems.
	This class will contain $\mathcal P \cap \mathcal K$, where
	$\mathcal P$ is the class of problems decidable in polynomial
	time and $\mathcal K$ is the class of problems definable without
	constants (See~\cite{KOIRAN97b} and Definition~\ref{K} below).
	Moreover:
\begin{theorem} \label{th1}
	The implications (a) $\Rightarrow$ (b)  $\Rightarrow$ (c) 
	 $\Rightarrow$ (d) are true:
	 \begin{itemize}
	 \item [(a)] The $\tau$-conjecture for polynomials.
	 \item [(b)] $\forall d$, $p_d$ is ultimately hard to compute.
	 \item [(c)] $\upnp$ over $\mathbb C$.
	 \item [(d)] $\pnp$ over $\mathbb C$.
	 \end{itemize}
\end{theorem}
\par
	The implication  (a) $\Rightarrow$ (b)  $\Rightarrow$ (d) appears
	in ~\cite{BCSS}, the hypothesis (c) in-between is new. It is at 
	least as likely as the $\tau$-conjecture, while still implying
	$\pnp$. 
\par
	We will also show a $\mathcal{NP}$-hardness result for 
	the class $\mathcal{UP}$: there is a 
	structured problem $\nullstellensatz \in \mathcal{NP} \cap \mathcal K$,
	such that:
\begin{theorem} \label{th2}
	$\upnp$ over $\mathbb C$ if and only if $\nullstellensatz \not
	\in \mathcal{UP}$ over  $\mathbb C$.
\end{theorem}
\par
	The problem $\nullstellensatz$ is precisely the (structured)
	Hilbert Nullstellensatz, known to be $\mathcal{NP}$-complete
	over $\mathbb C$~(\cite{BCSS}).
\medskip
\par
	This paper was written while the author was visiting Mathematical
	Sciences Research Institute in Berkeley. The author also wishes
	to thank Pascal Koiran and Steve Smale for their comments and
	suggestions.

\section{Background and Notations}

	Recall from $\cite{BCSS}$ that $\mathbb C^{\infty}$ is the
	disjoint union 
\[
	\mathbb C^{\infty} = \bigsqcup_{i=0, 1, \cdots} \mathbb C^i
\]
	This means that there is a well-defined {\em size} function, 
\[
\begin{array}{crcl}
	\size: & \mathbb C^{\infty} & \rightarrow & \mathbb N \\
	       & x & \mapsto & \size(x) = i \text{\ such that \ }
	       x \in \mathbb C^i
\end{array}
\]
\par
	A decision problem $X$ is a subset of $\mathbb C^{\infty}$.
	It is in the class $\mathcal P$ if and only if there is a
	machine $M$ over $\mathbb C$, that terminates for any input
	$x$ in time bounded by a polynomial on $\size(x)$, and 
	such that
\[
	M(x) = 0 
	\ \Leftrightarrow \ 
	x \in X
\]
	where $M(x)$ is the result of running $M$ with input $x$. Without
	loss of generality we may assume that $M(x) \in \{ 0;1 \}$.
\medskip
\par
	Under some circumstances, it is possible to assume that the
	machine $M$ above has only coefficients 0 or 1 (This is called
	a {\em constant-free} machine). However, one
	may have to replace problem $X$ over $\mathbb C$ by problem
	$X \cap \mathbb Z$ over $\mathbb Z$, with unit cost. (This is
	the contents of Propositions~3 and~9 of Chapter~7 of ~\cite{BCSS}).
	In order to avoid this technical complication and keep the
	same problem over $\mathbb C$, we will follow another approach
	to Elimination of Constants.
\par
	This approach was introduced by Koiran in ~\cite{KOIRAN97b}.
	The idea is to consider only machines for a subclass of problems. 
	This subclass will contain most of the interesting examples,
	while precluding pathological cases such as $X = \{ \pi \}$. 
\begin{definition}[Koiran] \label{K}
	A problem $L$ is said to be {\em definable without constants} if
	for each input size $n$ there is a formula $F_n$ in the first
	order theory of $\mathbb C$ such that $0$ and $1$ are the only
	constants occurring in $F_n$, and for any $x \in \mathbb C^n$,
	$x \in L$ if and only if $F_n(x)$ is true (there is no restriction
	on the size of $F_n$.
\end{definition}
\par
	For future reference, we quote below Theorem~2 of~\cite{KOIRAN97b}.
	The original statements of both Definition~\ref{K} and Theorem~\ref
	{ELIM} are actually more general (for any algebraically closed field
	of characteristic 0).
\begin{theorem}[Koiran] \label{ELIM}
	Let $L \subseteq K^\infty$ be a problem which is definable without
	constants. If $L \in \mathcal P$, $L$ can be recognized in
	polynomial time by a constant-free machine.
\end{theorem}
\par
	The class of all the problems  definable without constants
	will be denoted by $\mathcal K$.
\medskip
\par

	We will need crucially in the sequel the notion of a {\em
	structured problem}. A structured problem is a pair 
	$(X, X\yes)$, $X\yes \subseteq X \subseteq \mathbb C^{\infty}$.
	A non-structured problem $X$ can always be written as 
	the structured problem $(\mathbb C^{\infty}, X)$.
	The class $\mathcal{UP}$ will be meaningful only as a 
	class of structured problems. But first of all, recall that

\begin{definition} A structured problem $(X, X\yes)$ belongs to the
	class $\mathcal {P}$ if and only if $X \in \mathcal P$ and
	$X\yes \in \mathcal P$.
\end{definition}

\begin{definition} A structured problem $(X, X\yes)$ belongs to the
	class $\mathcal {K}$ if and only if $X \in \mathcal K$ and
	$X\yes \in \mathcal K$.
\end{definition}

\begin{definition} A structured problem $(X, X\yes)$ belongs to the
	class $\mathcal {NP}$ if and only if:
\begin{itemize}
	\item [(1)] The problem $X$ 
		belongs to the class $\mathcal P$.
	\item [(2)] There is a machine $M$ with input $x,g$ such that
\[
	x \in X \text{ and } \exists g \in \mathbb C^{\infty} 
	\text { s.t. } M(x,g)=0
	\
	\Leftrightarrow
	\ 
	x \in X\yes
\]
	\item [(3)] Furthermore, there is a polynomial $p$ such that,
	for all $x \in X\yes$, there is $g \in \mathbb C^{\infty}$ such
	that $M(x,g)=0$ and the running time of $M$ with input $x,g$ is
	no more than $p(\size(x))$.
\end{itemize}
\end{definition}

\begin{example}\label{ex1}
	Let $\hn$ be the class of all lists $(m, n, f_1, \cdots,f_m)$
	where $f_1, \cdots,$ $f_m$ are polynomials in $n$ variables. Each 
	polynomial $f = \sum f_I x^I$ is represented sparsely 
	by a list of monomials $(S, m_1, \cdots,
	m_S)$, where each monomial is a list $(f_I, I_1, \cdots, I_n)$.
\par
	An important convention to have in mind: integers appearing in
	the definition of a problem should be represented in bit
	representation. In this case, $m, n, S, I_j$ are all lists of
	zeros and ones. Complex values are represented by one complex
	number. With this convention, $\hn$ is clearly in the class
	$\mathcal P$.
\par
	We also define $\hn \yes$ as the subset of polynomial systems 
	in $\hn$ that have a common root over~$\mathbb C$.
\par	The definition above of the structured problem $\nullstellensatz$ 
	can be translated into first order constant-free formulae over
	$\mathbb C$. Therefore, $\nullstellensatz \in \mathcal K$.
	It is also $\mathcal{NP}$-complete over the complex numbers
	(Theorem~1 in Chapter~5 of ~\cite{BCSS}). 
\end{example}

\begin{example}\label{ex2}
	Let 
\begin{eqnarray*}
	X &=& \left\{ (m,x) \in \mathbb N \times \mathbb C \right\} \\
	X\yes &=& 
	      \left\{ (m,x) \in X \text{ such that } x \in \{ 1, 2, \cdots, m \}
	\right\} \\
\end{eqnarray*}
	with the convention that $m$ is in bit representation, while $x$ is
	a complex number. Hence, $\size((m,x)) = O(1 + \lceil \log_2 (m)
	\rceil)$. Then the problem $(X, X\yes)$ is in $\mathcal {NP}$ over
	$\mathbb C$. The machine $M(x,g)$ can be constructing by
	guessing the bit decomposition $g_i$ of $x$, and computing
	$x - \sum g_i 2^i$.
\par
	Again, $(X,X\yes)$ is definable without constants.
\end{example}

\section{Construction of the class $\mathcal {UP}$}
\par
	In Chapter 7 of~\cite{BCSS}, it is proved that if the problem
	$(X, X\yes)$ from Example~\ref{ex2} would happen to
	belong to the class $\mathcal {P}$, then condition (b) in Theorem
	~\ref{th1} would be false. Therefore (b) implies \pnp over 
	$\mathbb C$.
\medskip
\par
	The class $\mathcal{UP}$ will be constructed by abstracting the
	same reasoning. The construction relies on some geometric
	properties of structured problems in $\mathcal P$. The notation
	that follows will be used in the sequel:
\par
	Let $(X, X\yes)$ be a structured problem  with $X \in \mathcal P$.
	We denote by $X \cap \mathbb C^i$ the set  
	$\{ x \in X: \size(x) = i \}$  
	of size $i$ instances of the problem. Then we write
	$\overline {X \cap \mathbb C^i}$ for its Zariski closure 
	over $\mathbb C$. We can define a new object associated to 
	$X$ as:
\[
	\overline{X} = \bigsqcup_{i=0,1, \cdots} \overline {X \cap \mathbb C^i}
\]
\par
	We can think of $\overline{X}$ as the {\em closure}
	of $X$, indeed it is the smallest `closed' problem containing 
	$X$. Remark that in Examples~\ref{ex1} and~\ref{ex2}, we have
	respectively $X=\overline{X}$ and $\hn=\overline{\hn}$. 
\par
	We can also decompose each Zariski-closed set $\overline {X \cap
	\mathbb C^i}$ into a finite union of irreducible components
	(affine varieties). Thus it makes sense to write $\overline{X}$ as
	the countable union:
\[
	\overline{X} = \bigcup X_j
\]
	where each $X_j$ is an affine variety lying in some $\mathbb C^s$,
	where $s = \size (x), x \in X_j$. We can further define:
\begin{eqnarray*}
	X_j\yes &=& X_j \cap X\yes \\
	X_j\no &=& X_j \setminus X\yes 
\end{eqnarray*}
\par
\begin{figure}
\centerline{
\resizebox{4cm}{7cm}{\includegraphics{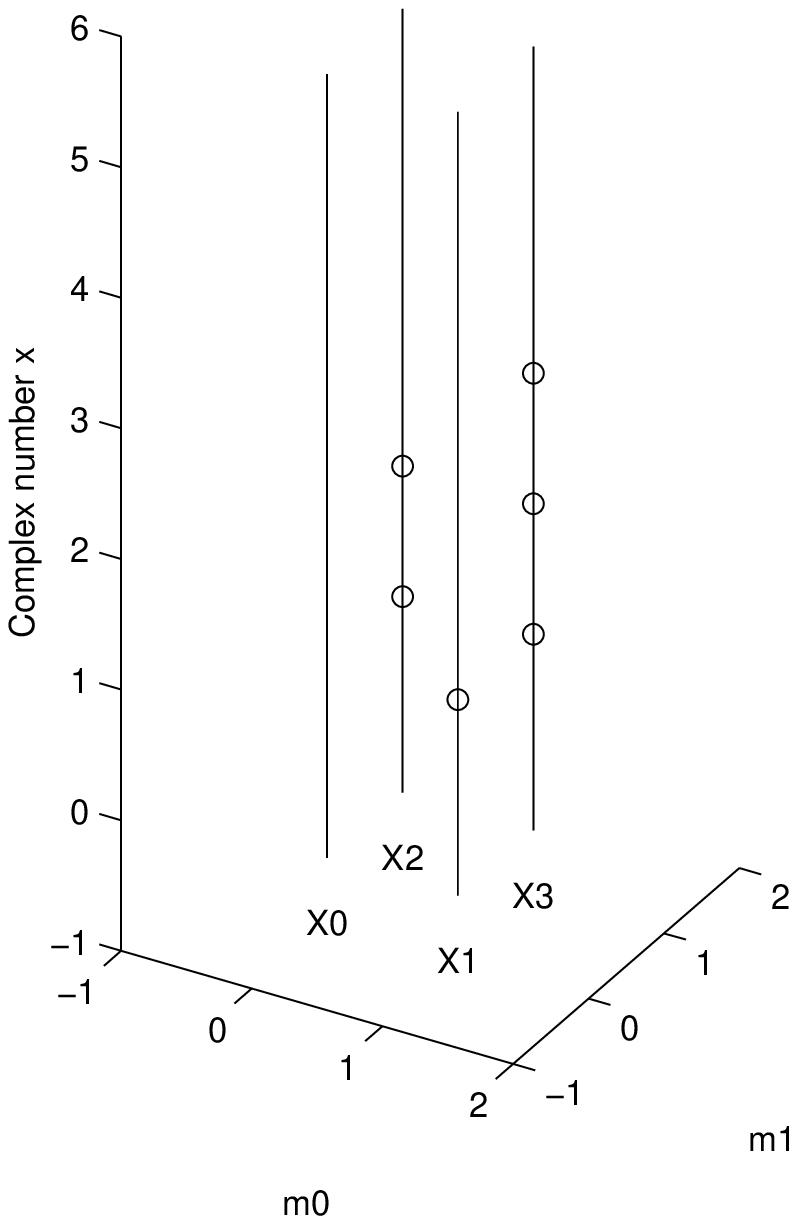}}
\hspace{2cm}
\begin{minipage}[b]{4cm}
{\sf \footnotesize This is Problem $(X,X\yes)$
  from Example~\ref{ex2}, restricted to the inputs
  $(m_0, m_1, x)$ of size $3$. 
  $X$ is represented by the four (complex !) lines and 
  $X\yes$ by the dots. 
  Each of the complex lines is irreducible, and hence
  corresponds to a different $X_i$.}
  \vspace{1cm}
\end{minipage}
\hspace{2cm}
}
\caption{\label{fig1} $(X,X\yes)$ from Example~\ref{ex2}}
\end{figure}

\par
	(See Figure~\ref{fig1}).  Using this notation, 
\begin{definition}\label{dup}
	The class $\mathcal{UP}$ is the class of all structured problems
	$(X, X\yes)$ such that $X \in \mathcal{P}$ and for all $X_i$,
	there is a non-zero polynomial $f_i \in \mathbb Z[x_1, \cdots,
	x_{s_i}]$, where $s_i = \size (x)$ for $x \in X_i$, with the following
	properties:
\begin{itemize}
	\item [(1)] $\tau(f_i)$ is polynomially bounded in $S_i$.
	\item [(2)] $X_i\yes \subseteq Z(f)$ or $X_i\no \subseteq Z(f)$
\end{itemize}
\end{definition}
\medskip
\par

\begin{proposition}\label{pup}
	$\mathcal P \cap \mathcal K \subseteq \mathcal {UP}$
\end{proposition}

\begin{proof}[Proof of Proposition~\ref{pup}]
	Let $(X,X\yes)$ be in $\mathcal P \cap \mathcal K$. 
	Let $M = M(x)$ be the
	machine that recognizes $x \in X\yes$ in polynomial time,
	where the input $x$ is assumed to be in $\overline{X}$.
	Although it is possible that an $x \in X_i$ is not in $X$,
	it is still possible to recognize $x \in X\yes$ in polynomial
	time. Indeed, $X$ is also in $\mathcal P$. The
	machine $M(x)$ will check $x \in X$
	and $x \in X\yes$. 
\par
	Now we apply elimination of constants (Theorem~\ref{ELIM}),
	and choose $M$ to be constant-free. 	
\par
	The nodes of the machine $M$ are supposed to be numbered.
	Given an input $x$, the {\em path} followed by input $x$
	is the list of nodes traversed during the computation of
	$M(x)$.
\par
	When the input is restricted to one of the affine varieties
	$X_i$'s, we can define the canonical path (associated to
	$X_i$ as the path followed
	by the generic point of $X_i$. This corresponds to the following
	procedure:
\par
	At each decision node, at time $T$, branch depends upon an 
	equality $F^T(x) = 0$, where $x$ is the original input. The
	polynomial $F$ can be computed within the machine running time.
	In case $F^T(x) = 0$ for all $x \in X_i$, we follow the
	Yes-path and say that this branching is trivial.
\par
	If not, we follow the no-path and say that this branching is
	non-trivial. The fact that $X_i$ is a variety is essential
	here, since it guarantees that only a codimension $\ge 1$ 
	subset of inputs may eventually follow the Yes-path at this
	time.
\par
	The set of inputs that do NOT follow the canonical path can
	be described as the zero-set of
\[
	f_i = \prod F^T
\]
	where the product ranges over the non-trivial branches only.
	The polynomial $f_i$ can be computed in at most twice the 
	running time of the machine $M$ restricted to $X_i$. By
	hypothesis, this is polynomial time in the size of $x \in X_i$. 
\par
	Since we assumed that $M$ returns only $0$ or $1$, the set
	of the inputs that follow the canonical path (i.e. $Z(f_i)$)
	is either
	all in $X_i\yes$ or all in its complementary $X_i\no$.
\par
	There are now two possibilities. First possibility, $X_i\yes$ 
	has measure zero in $X_i$, and therefore it must be contained
	in $Z(f_i)$. Second possibility, $X_i\yes$ has non-zero
	measure, hence it contains the complementary of $Z(f_i)$, and
	hence $X_i\no$ is a subset of $Z(f_i)$.
\end{proof}

\section {Proof of the Theorems}

\begin{proof}[Proof of Theorem~\ref{th1}]
\medskip \ {}
\par
(a) $\Rightarrow$ (b) is trivial, refer to ~\cite{BCSS} Chapter 7.
\medskip
\par
(b) $\Rightarrow$ (c): Let $(X,X\yes)$ be the problem in Example~\ref{ex2}.
	Since $X_i\no$ is generic in $X_i$, all inputs in $X_i\yes$ should
	escape the canonical path. Hence, if $f_d$ is the polynomial that
	defines the canonical path, $f_d(i) = 0$ for $i=1,2, \cdots, d$.
	But then it cannot be evaluated in time polylog($d$), by hypothesis
	(b). Hence, under the assumption (b), the problem $(X,X\yes)$ is
	not in $\mathcal {UP}$. It does belong to $\mathcal {NP} \cap 
	\mathcal K$, so $\upnp$.
\medskip
\par
(c) $\Rightarrow$ (d) : Using Theorem~\ref{th2}, Condition (c) implies
	that $\nullstellensatz \not \in \mathcal{UP}$. However, since
	$\nullstellensatz \in \mathcal K$, Proposition~\ref{pup}
	implies $\nullstellensatz \not \in \mathcal P$. Hence \pnp
	over $\mathbb C$.
\end{proof}

\begin{proof}[Proof of Theorem~\ref{th2}]

	Let $(X, X\yes) \in \mathcal{NP} \cap \mathcal{K}$ and assume that
	$\nullstellensatz \in \mathcal{UP}$. We have to show 
	that $(X, X\yes) \in \mathcal{UP}$. 
\par
	For each $X_i$, one can embed $(X_i, X_i\yes)$ into some
	$(\hn_i, \hn_i\yes)$ as follows:
\par
	Let $M = M(x)$ be the deterministic polynomial time
	machine to recognize $X$, and let $N=N(x,g)$ be the
	non-deterministic polynomial time machine to recognize
	$X\yes$. We can assume without loss of generality that
	$M$ and $N$ are constant-free (Theorem~\ref{ELIM}).
\par
	Let $T$ be the maximum running time of $M$ and $N$
	when the input is
	restricted to $X_i$. Let $\phi(x)$ be the combined Register Equations
	of machines $M$ and $N$ for time $T$ 
	(Theorem~2 in Chapter~3 of~\cite{BCSS}). 
	Thus, $\phi(x)$ is a system of
	polynomial equations with integer coefficients and indeterminate
	coefficients $x_1, x_2, \cdots$.
	The polynomial system $\phi(x)$ can be constructed in polynomial
	time from $x$, and the size
	of $\phi(x)$ is polynomially bounded by the size of $x$.
\medskip
\par
	We claim that $\phi(X_i)$ is contained in some $\hn_j$, and
	that in that case $\phi(X_i\yes) \subseteq \hn_j\yes$ and
	$\phi(X_i\no) \subseteq \hn_j\no$.
\par
	Indeed, $X_i \subseteq \mathbb C^s$ for some $s$, and
	$\phi(\mathbb C^s) \subseteq HN_j$ for some $j$.
	Then $x \in X_i$ belongs to $X\yes$ if and only if the corresponding
	$\phi(x)$ has a solution over $\mathbb C$.
\medskip
\par
	We now distinguish two cases:
\par	Case 1: $\hn_j\yes$ has measure zero in $\hn_j$. Thus
	$\hn_j\yes \subseteq Z(\hat f_j)$ for an easy-to-compute 
	polynomial $\hat f_j$. In that case,
	since $X_i\yes$ gets mapped into $\hn_j\yes$, the composition
	$f_i = \hat f_j \circ \phi$ gives the polynomial associated to $X_i$.
\par	Case 2: $\hn_j\no$ has measure zero in $\hn_j$. Thus
        $\hn_j\no \subseteq Z(\hat f_j)$ for an easy-to-compute 
	polynomial $\hat f_j$. In that case,
        since $X_i\no$ gets mapped into $\hn_j\no$, 
        $f_i = \hat f_j \circ \phi$ is the polynomial associated to $X_i$.
\end{proof}

\section{Ultimate Complexity}

	Let $(Y,Y\yes)$ be a problem over $\mathbb C$, definable without
	constants and with $Y$ semi-decidable (i.e. $Y$ is the halting
	set of some machine). The closure $\overline Y$ is well-defined and
	can be written as a countable union of irreducible varieties
	$Y_i$. 
\par
	For any machine $M$ to solve $(Y,Y\yes)$, one can produce
	a family of polynomials $f_i$, vanishing on the set of inputs
	that follow the canonical-path of $M$ restricted to $Y_i$.
	As in item (2) of Definition~\ref{dup}, we have
\[
	Y_i\yes \subseteq Z(f_i) \text{\ or\ } Y_i\no \subseteq Z(f_i)
\]
\par
	Also, for each input size $s$, one has a finite number
	of indices $i$ corresponding to components i
	$Y_i \subseteq \overline Y$ of  
	size-$s$ input. We can thus maximize over those
	indices $i$:
\[
	u_M (s) = \max_{i:Y_i \subseteq \mathbb C^{i}} \tau(f_i)
\]
\par
	This invariant may be called `ultimate running time', and
	is a lower bound (up to a constant) for the worst-case
	running time of $M$. As with ordinary complexity theory,
	one can define the `ultimate complexity' class of a problem
	as the class of functions $u: \mathbb N \rightarrow \mathbb R$
	such that $\exists M, c>0: \forall x u_M(x) \le c u(x)$ and 
	$M$ recognizes $(Y, Y\yes)$. This provides notions such as
	`ultimate logarithmic time' or `ultimate exponential time'.
\par
	In~\cite{M99}, a similar construction is used to obtain
	lower bounds for some specific decision problems. Those
	problems, however, had a very simple geometric structure
	(for each `input size', $X\yes$ was a finite set in 
	$\mathbb C$). The motivation of this paper was to extend
	some of the ideas therein and in Chapter~7 of~\cite{BCSS}
	to non-codimension-1 problems.

\bibliographystyle{plain}
\bibliography{up}

\end{document}